\documentclass[11pt]{article}
\usepackage{amssymb,amsmath,amsfonts,mathtools}
\usepackage{epsfig,graphicx}
\usepackage{latexsym}
\newtheorem{theorem}{Theorem}
\newtheorem{lemma}[theorem]{Lemma}
\newtheorem{coro}[theorem]{Corollary}

\def\qed
 {\ifhmode\unskip\nobreak\hfill$\Box$\medskip\fi
 \ifmmode\eqno{\Box}\fi}

\def\ex{{\rm ex}}
\def\EE{{\mathcal E}}

\def\ed{{\rm ed}}

\def\LL{{\mathcal L}}

\def\PP{{\mathcal P}}

\begin{document}

\title{{
A proof of the stability of extremal graphs,\\
 Simonovits' stability from Szemer\'edi's regularity}}

\author{Zolt\'an F\"uredi
\thanks{
Research supported in part by the Hungarian National Science Foundation OTKA 104343,
 by the Simons Foundation Collaboration Grant \#317487,
and by the European Research Council Advanced Investigators Grant 267195.
\newline\indent
{\it 2010 Mathematics Subject Classifications:}
05C35. \hfill \jobname
\newline\indent
{\it Key Words}:  Tur\'an number,  extremal graphs, stability.  \hfill\today
}
}
\date{Alfr\'ed R\'enyi Institute of Mathematics,\\ 13--15 Re\'altanoda Street, 1053 Budapest, Hungary. \\
E-mail: {\tt z-furedi@illinois.edu}}


\maketitle

\maketitle

\begin{abstract}
The following sharpening of Tur\'an's theorem is proved.
Let $T_{n,p}$ denote the complete $p$--partite graph of order $n$ having the maximum number of edges.
If $G$ is an $n$-vertex $K_{p+1}$-free graph with
  $e(T_{n,p})-t$ edges then there exists an (at most) $p$-chromatic subgraph $H_0$ such that
  $e(H_0)\geq e(G)-t$.

Using this result we present a concise, contemporary proof (i.e., one using
Szemer\'edi's regularity lemma) for the classical stability result
of Simonovits~\cite{Sim_Tihany}.

\end{abstract}

\section{The Tur\'an problem} 

Given a graph $G$ with vertex set $V(G)$ and vertex set $\EE(G)$ its number of edges is denoted by $e(G)$.
The neighborhood of a vertex $x\in V$ is denoted by $N(x)$, note that $x\notin N(x)$.
For any $A\subset V$ the restricted neighborhood $N_G(x|A)$ stands for $N(x)\cap A$.
Similarly, $\deg_G(x|A):=|N(x)\cap A|$.
If the graph is well understood from the text we leave out subscripts.
The {\it Tur\'an graph} $T_{n,p}$ is the largest $p$-chromatic graph having $n$ vertices, $n,p\geq 1$.
Given a partition $(V_1, \dots, V_p)$ of $V$ the {\it complete multipartite} graph
 $K(V_1, \dots, V_p)$ has vertex set $V$ and all the edges joining distinct partite sets.
$A\bigtriangleup B$ stands for the symmetric difference of the sets $A$ and $B$.
For further notations and notions undefined here see, e.g., the monograph of Bollob\'as~\cite{BolExt}.

Tur\'an~\cite{TuranA} proved that if an $n$ vertex graph $G$ has at least
 $e(T_{n,p})$ edges then it contains a complete subgraph $K_{p+1}$, except if $G=T_{n,p}$.
Given a class of graphs $\LL$, a graph $G$ is called $\LL$-{\it free} if it does not contain any subgraph isomorphic to any member of $\LL$.
The {\it Tur\'an number} $\ex(n,\LL)$ is defined as the largest size of an $n$-vertex, $\LL$-free graph.
Erd\H os and Simonovits~\cite{ErdSimLim} gave a general asymptotic for the Tur\'an number as follows.
Let $p+1:=\min \{\chi(L): L\in\LL\}$. Then
\begin{equation}\label{eq1}
\ex(n,\LL)=\left(1-\frac{1}{p}\right)\binom{n}{2}+o(n^2) \quad{\rm as} \quad n\to \infty.
  \end{equation}
They also showed that if $G$ is an extremal graph, i.e., $e(G)=\ex(n, \LL)$,
 then it can be obtained from $T_{n,p}$ by adding and deleting at most $o(n^2)$ edges.
This result is usually called Erd\H os--Stone--Simonovits theorem, although it was proved first in~\cite{ErdSimLim},
 but indeed (\ref{eq1}) easily follows from a result of Erd\H os and Stone~\cite{ErdStone}.

The aim of this paper is to present a new proof for the following stronger version of (\ref{eq1}), a structural stability theorem,
 originally proved by Erd\H os and  Simonovits~\cite{ErdSimLim}, Erd\H os~\cite{Erd_Rome,Erd_Tihany}, and Simonovits~\cite{Sim_Tihany}.
For every $\varepsilon>0$ and forbidden subgraph class $\LL$ there is a $\delta>0$, and $n_0$ such that
 if  $n>n_0$ and $G$ is an $n$-vertex, $\LL$-free graph with
\begin{equation*}
 e(G)\ge\left(1-\frac{1}{p}\right)\binom{n}{2}-\delta n^2,
   \end{equation*}
then
\begin{equation}\label{eq2}
|\EE(G_n) \bigtriangleup \EE(T_{n,p})|\le \varepsilon n^2.
 \end{equation}
I.e., one can change (add and delete) at most $\varepsilon n^2$ edges of $G$
 and obtain a complete $p$-partite graph.
In other words, if an $n$-vertex $\LL$-free graph $G$ is almost extremal,
$\min \{ \chi(L): L\in \LL\}=p+1$, then the structure of $G$ is close to a $p$-partite Tur\'an graph.
This result is usually called Simonovits' stability of the extremum.

Our main tool is a very simple proof for the case $\LL=\{K_{p+1}\}$.

Stability results are usually  more important than their extremal counterparts.
That is why there are so many investigations concerning the {\it edit distance} of graphs.
Let $G_1=(V,\EE_1)$ and $G_2=(V,\EE_2)$ be two (finite, undirected) graphs on the same vertex set.
The {\it edit distance} from $G_1$ to $G_2$ is $\ed(G_1,G_2):=|\EE_1 \bigtriangleup  \EE_2|$.
Let $\PP$ denote a class of graphs and $G$ be a fixed graph.
The edit distance from $G$ to $\PP$ is $\ed(G,\PP)=\min\{\ed(G,F):F\in \PP,V(G)=V(F)\}$.
This notion was explicitly introduced in~\cite{AxeMartin2006}, Alon and Stav~\cite{AloSta} proved connections with Tur\'an theory.
For more recent results see Martin~\cite{Martin2015}.

\section{How to make a $K_{p+1}$-free graph $p$-chromatic}

Ever since Erd\H os~\cite{ErdBipartite} observed that one can always delete at most $e/2$ edges from any graph $G$
 to make it bipartite there are many generalizations and applications of this
(see, e.g., Alon~\cite{Alo} for a more precise form).
Here we prove a version dealing with a narrower class of graphs.
Recall that $e(T_{n,p}):= \max\{ e(K(V_1, \dots, V_p)): \sum |V_i|=n\}$, the maximum size of a $p$-chromatic graph.

\begin{theorem}\label{th1}
Suppose that $K_{p+1}\not\subset G$, $|V(G)|=n$, $t\geq 0$, and
 $$ e(G)= e(T_{n,p})-t.$$
Then there exists an (at most) $p$-chromatic subgraph $H_0$,
 $\EE(H_0)\subset \EE(G)$ such that
$$  e(H_0)\geq e(G)-t. $$
\end{theorem}

\begin{coro}[Stability of $\ex(n, K_{p+1})$]\label{co1}\enskip
Suppose that $G$ is $K_{p+1}$-free with $e(G)\geq e(T_{n,p})- t$.
Then there is a complete $p$-chromatic graph $K:=K(V_1, \dots, V_p)$ with $V(K)=V(G)$,
  such that
$$   |\EE(G)\bigtriangleup \EE(K)|\leq 3t. $$
\end{coro}
Indeed, delete $t$ edges of $G$ to obtain the $p$-chromatic $H_0$.
Since $e(H_0)\geq e(T_{n,p})-2t$ one can add at most $2t$ edges to make it a complete $p$-partite graph.
(Here $V_i=\emptyset$ is allowed). \qed

There are other more exact stability results, e.g.,
Hanson and Toft~\cite{HanTof} showed that for  $t< n/(2p)-O(1)$ the graph $G$ itself is $p$-chromatic,
 there is no need to delete any edge.
Some results of E. Gy\H ori~\cite{Gyo} implies a stronger form, namely that $e(H_0)\geq e(G)-O(t^2/n^2)$.
Erd\H os, Gy\H ori, and  Simonovits~\cite{ErdGyoSim} considers only dense triangle-free graphs.
The advantage of our Theorem~\ref{th1} is that it contains no $\varepsilon, \delta, n_0$,
 it is true for every $n$, $p$ and $t$.

The inequality in Corollary~\ref{co1} is simple because we estimate the edit distance of $G$ from
 a not necessarily balanced $p$ partite graph $K$.
If we are interested in $\ed(G,T_{n,p})$ then we can use the following inequality.
If $e(K((V_1, \dots, V_p))\geq e(T_{n,p})-2t$, then a simple calculation shows that
the sizes of $V_i$'s should be 'close' to $n/p$ (more exactly we get
$4t\geq \sum_i (|V_i|-(n/p))^2$) and hence
\begin{equation}\label{eq:co2}
  \ed(K,T_{n,p})\leq n \sqrt{t/p}
\end{equation}

\noindent
{\bf Proof of Theorem~\ref{th1}.} \quad
We find the large $p$-partite subgraph $H_0\subset G$ by analyzing
Erd\H os' degree majorization algorithm~\cite{ErdTur} what he used to prove Tur\'an's theorem.
Our input is the $K_{p+1}$-free graph $G$ and the output is a partition $V_1, V_2, \dots, V_p$
 of $V(G)$ such that  $\sum_i e(G|V_i) \leq t$.

Let $x_1\in V(G)$ be a vertex of maximum degree and  let $V_1:= V\setminus N(x_1)$, $V_1^+:= V\setminus V_1$.
Note that $x_1\in V_1$ and $\deg(x)\leq |V_1^+|$ for all $x\in V_1$. Hence
\begin{equation*}
   2e(G|V_1)+ e(V_1, V_1^+ )=\sum_{x\in V_1} \deg(x)\leq |V_1||V_1^+|.
   \end{equation*}
In general, define $V_0^+:=V(G)$ and let $x_i$ be a vertex of maximum degree of the graph  $G|V_{i-1}^+$,
 let $V_i:= V_{i-1}^+\setminus N(x_i)$, $V_i^+:=V(G)\setminus (V_1\cup \dots \cup V_i)$.
We have $x_i\in V_i$, $\deg(x_i, V_{i-1}^+)=|V_i^+|$ and
\begin{equation}\label{eq4}
   2e(G|V_i)+ e(V_i, V_i^+ )=\sum_{x\in V_i} \deg(x| V_{i-1}^+)\leq |V_i||V_i^+|.
   \end{equation}
The procedure stops in $s$ steps when no more vertices left, i.e., if $V_1\cup \dots \cup V_s=V(G)$.
Note that $s\leq p$ because $\{ x_1, x_2, \dots , x_s\}$
 span a complete graph.

Add up the left hand sides of (\ref{eq4}) for $1\leq i\leq s$, we get  $e(G)+ \left(\sum_i e(G|V_i)\right)$.
The sum of the right hand sides is exactly  $e(K(V_1, V_2, \dots , V_s))$.
We obtain
  $$
  e(T_{n,p})-t+ \left(\sum_i e(G|V_i)\right)= e(G)+ \left(\sum_i e(G|V_i)\right)\leq e(K(V_1, V_2, \dots , V_p))
    \leq e(T_{n,p})
    $$
implying  $\sum_i e(G|V_i)\leq t $.  \qed

\section{Az application of the Removal Lemma}

We only need a simple consequence of Szemer\'edi's Regularity Lemma.
Recall that the graph $H$ contains a homomorphic image of $F$ if there is a mapping
 $\varphi:V(F)\to V(H)$ such that the image of each $F$-edge is an $H$-edge.
There is a $\varphi:V(F)\to V(K_s)$ homomorphism if and only if $s\geq \chi(H)$.
If there is no any $\varphi:V(F)\to V(H)$ homomorphism then $H$ is called  ${\rm hom}(F)$-free.

\begin{lemma}[A simple form of the Removal Lemma]\label{le1}\enskip
For every  $\alpha>0$ and graph $F$ there is an $n_1$ such that
 if  $n>n_1$ and $G$ is an $n$-vertex, $F$-free graph then it contains a ${\rm hom}(F)$-free subgraph $H$
 with $e(H)> e(G)-\alpha n^2$.
   \end{lemma}

This means that $H$ does not contain any homomorphic image of $F$ as a subgraph, especially
 if $\chi(F)=p+1$ then $H$ is $K_{p+1}$-free.
The Removal Lemma can be attributed to Ruzsa and Szemer\'edi~\cite{RuzSze}.
It appears in a more explicit form in \cite{ErdFraRod} and \cite{Fur}.
For a survey of applications of Szemer\'edi's regularity lemma in graph theory see Koml\'os-Simonovits~\cite{KomlosSim}
 or Koml\'os-Shokoufandeh-Simonovits-Szemer\'edi~\cite{KomShoSimSze}.

\noindent
{\bf Proof of (\ref{eq2})} using Lemma~\ref{le1} and Corollary~\ref{co1}.\enskip
Suppose that $F\in \LL$, $\chi(F)=p+1$ and $\alpha>0 $ an arbitrary real.
Suppose that $G$ is $F$-free with $n> n_1(F, \alpha)$ and $e(G)> e(T_{n,p})-
 \alpha n^2$.
We have to show that the edit distance of $G$ to $T_{n,p}$ is small.
First we claim that the edit distance of $G$ to a complete $p$--partite graph $K(V_1, \dots, V_p)$ is at most $7\alpha n^2$.
Indeed, using the Removal Lemma we obtain a $K_{p+1}$-free subgraph $H$ of $G$ such that
  $e(H)>e(G)-\alpha n^2 > e(T_{n,p})- 2\alpha n^2$.
Apply Theorem~\ref{th1} to $H$ we get a $p$--partite $H_0$ with $e(H_0)> e(T_{n,p})- 4\alpha n^2$.
Then Corollary~\ref{co1} yields a $K:=K(V_1, \dots, V_p)$ with $\ed(K,H) < 6\alpha n^2$, giving
 $\ed(K, G)\leq 7\alpha n^2$.

Since  $e(K)\geq e(H_0)> e(T_{n,p})-4\alpha n^2$, we can use (\ref{eq:co2}) with $t=2\alpha n^2$
 to get $\ed(K, T_{n,p})\leq n^2 \sqrt{2\alpha/p}$.
This completes the proof that $\ed(G,T_{n,p})\leq(7\alpha + \sqrt{2\alpha/p})n^2$.
\qed

\noindent
{\bf Acknowledgments.}\quad
The author is greatly thankful to 
 M.~Simonovits for helpful conversations.
\newline
This result was first presented in a public lecture at Charles University, Prague, July 2006.
Since then there were several references to it, e.g., in~\cite{Mub}.


\end{document}